# Generalized inverses of Markovian kernels in terms of properties of the Markov chain

Jeffrey J. Hunter

*School of Computing and Mathematical Sciences, Auckland University of Technology, Private Bag 92006, Auckland 1142, New Zealand*

## Abstract

All one-condition generalized inverses of the Markovian kernel $I - P$, where $P$ is the transition matrix of a finite irreducible Markov chain, can be uniquely specified in terms of the stationary probabilities and the mean first passage times of the underlying Markov chain. Special sub-families include the group inverse of $I - P$, Kemeny and Snell's fundamental matrix of the Markov chain and the Moore-Penrose g-inverse. The elements of some sub-families of the generalized inverses can also be re-expressed involving the second moments of the recurrence time variables. Some applications to Kemeny's constant and perturbations of Markov chains are also considered.

## 1.     Introduction.

Let $P = [p_{ij}]$ be the transition matrix of a finite irreducible, discrete time Markov chain $\{X_n\}$, $(n \geq 0)$, with state space $S = \{1, 2,..., m\}$. Such chains have a unique stationary distribution $\{\pi_j\}$, $(1 \leq j \leq m)$. Let $T_{ij} = \min[n \geq 1, X_n = j \mid X_0 = i]$ be the first passage time from state $i$ to state $j$ (first return when $i = j$) and define $m_{ij} = E[T_{ij} \mid X_0 = i]$ as the mean first passage time from state $i$ to state $j$ (or mean recurrence time of state $i$ when $i = j$). It is well known that for finite irreducible chains all the $m_{ij}$ are well defined and finite.

Generalized matrix inverses (g-inverses) of $I - P$ are typically used to solve systems of linear equations and a variety of the properties of the Markov chain, in particular the $\{\pi_j\}$ and the $\{m_{ij}\}$ can be found in terms of g-inverses, either in matrix form or in terms of the elements of the g-inverse. What is not known is that the elements of every g-inverse of $I - P$ can be expressed in terms of the stationary probabilities $\{\pi_j\}$ and the mean first passage times $\{m_{ij}\}$ of the associated Markov chain. The key thrust to this paper is to first identify the parameters that characterise different sub-families of g-inverses of $I - P$. Then to assign to each sub family, so characterised, explicit expressions for the elements of the g-inverses in terms of the $\{\pi_j\}$ and the $\{m_{ij}\}$.

## 2.   Generalized inverses of the Markovian kernel $I - P$.

A g-inverse of a matrix $A$ is any matrix $A^-$ such that $AA^-A = A$. Matrices $A^-$ are called "one condition" g-inverses or "equation solving" g-inverses because of their use in solving systems of linear equations.

If $A$ is non-singular then $A^-$ is $A^{-1}$, the inverse of $A$, and is unique. If $A$ is singular, $A^-$ is not unique. Typically, in the equations that we wish to solve, we only need a one-condition g-inverse with the non-uniqueness being eliminated by the imposition of boundary conditions (such as $\sum_{i=1}^{m} \pi_i = 1$ in the case of finding stationary distributions, and $m_{ii} = 1/\pi_i$ in the case of mean first passage times).

The following theorem, due to Hunter, (1982), gives a procedure for finding all one-condition g-inverses of $I - P$.

**Theorem 1**: *Let P be the transition matrix of a finite irreducible Markov chain with m states and stationary probability vector $\pmb{\pi}^T = (\pi_1, \pi_2, \ldots, \pi_m)$.*
*Let $\pmb{e}^T = (1, 1, \ldots, 1)$ and $\pmb{t}$ and $\pmb{u}$ be any vectors.*
(a)   $I - P + \pmb{t}\pmb{u}^T$ *is non-singular if and only if $\pmb{\pi}^T \pmb{t} \neq 0$ and $\pmb{u}^T \pmb{e} \neq 0$.*
(b)   *If $\pmb{\pi}^T \pmb{t} \neq 0$ and $\pmb{u}^T \pmb{e} \neq 0$ then $[I - P + \pmb{t}\pmb{u}^T]^{-1}$ is a one condition g-inverse of $I - P$ and, further, all "one condition" g-inverses of $I - P$ can be expressed as*
   $A^{(1)} = [I - P + \pmb{t}\pmb{u}^T]^{-1} + \pmb{e}\pmb{f}^T + \pmb{g}\pmb{\pi}^T$ *for arbitrary vectors $\pmb{f}$ and $\pmb{g}$.*

□

Useful by-products of the proof of the above theorem were the following results:

$$[I - P + \pmb{t}\pmb{u}^T]^{-1}\pmb{t} = \frac{\pmb{e}}{\pmb{u}^T \pmb{e}}. \tag{2.1}$$

$$\pmb{u}^T [I - P + \pmb{t}\pmb{u}^T]^{-1} = \frac{\pmb{\pi}^T}{\pmb{\pi}^T \pmb{t}}. \tag{2.2}$$

Special multi-condition g-inverses of $A$ can also be considered by imposing additional conditions. Consider real conformable matrices $X$ (which in our context we assume to be square) such that:

    (Condition 1)    $AXA = A$.
    (Condition 2)    $XAX = X$.
    (Condition 3)    $(AX)^T = AX$.
    (Condition 4)    $(XA)^T = XA$.
    (Condition 5)    $AX = XA$.

Let $A^{(i,j,...,l)}$ be any matrix that satisfies conditions *(i), (j), ..., (l)* of the above itemised conditions. $A^{(i,j,...,l)}$ is called an - *(i, j, ..., l)* g-inverse of *A*, under the assumption that condition (1) is always included. Let $A\{i, j, ..., l\}$ be the class of all *(i, j, ..., l)* g-inverses of *A*.

The classification of g-inverses of the Markovian kernel $I - P$, can be done simply by means of the following results of Hunter, (1990).

**Theorem 2**: *If G is any g-inverse of $I - P$, where P is the transition matrix of a finite irreducible Markov chain with stationary probability vector $\pi^T$, then G can be uniquely expressed in parametric form as*

$$G \equiv G(\alpha, \beta, \gamma) = [I - P + \alpha\beta^T]^{-1} + \gamma e\pi^T, \tag{2.3}$$

*where $\alpha, \beta$, and $\gamma$ involve $2m - 1$ independent parameters with the property that*

$$\pi^T\alpha = 1 \text{ and } \beta^T e = 1. \tag{2.4}$$

*If $A \equiv I - (I - P)G$ and $B \equiv I - G(I - P)$ then*

$$A = \alpha\pi^T \text{ and } B = e\beta^T, \tag{2.5}$$

*so that* $\quad\quad\quad\quad\quad\quad \alpha = Ae \text{ and } \beta^T = \pi^T B. \tag{2.6}$

*Further, from (2.1) and (2.2), $G\alpha = (\gamma + 1)e$ and $\beta^T G = (\gamma + 1)\pi^T$,* $\tag{2.7}$

*so that* $\quad\quad\quad\quad\quad\quad \gamma + 1 = \pi^T G\alpha = \beta^T Ge = \beta^T G\alpha. \tag{2.8}$

*Also*

$$G \in A\{1, 2\} \Leftrightarrow \gamma = -1, \tag{2.9}$$
$$G \in A\{1, 3\} \Leftrightarrow \alpha = \pi / \pi^T\pi, \tag{2.10}$$
$$G \in A\{1, 4\} \Leftrightarrow \beta = e / e^T e = e / m, \tag{2.11}$$
$$G \in A\{1, 5\} \Leftrightarrow \alpha = e, \beta = \pi. \tag{2.12}$$

□

We further subdivide the classification given by (2.12) and define
$$G \in A\{1, 5a\} \Leftrightarrow \alpha = e, \tag{2.13}$$
$$G \in A\{1, 5b\} \Leftrightarrow \beta = \pi, \tag{2.14}$$
so that $G \in A\{1, 5\} \Leftrightarrow G \in A\{1, 5a\}$ and $G \in A\{1, 5b\}$.

**Theorem 3**: *Let G be any g-inverse of $I - P$, where P is the transition matrix of a finite irreducible Markov chain with stationary probability vector $\pi^T$. Suppose G has parametric form $G = G(\alpha, \beta, \gamma)$, where the parameters $\alpha, \beta$, and $\gamma$ are as given in Theorem 2.*
(a) $G \in A\{1, 5a\} \Leftrightarrow \alpha = e \Leftrightarrow Ge = ge$ *for some g.*
  *Further, if $Ge = ge$ for some g then $g = 1 + \gamma$.*
(b) $G \in A\{1, 5b\} \Leftrightarrow \beta = \pi \Leftrightarrow \pi^T G = h\pi^T$ *for some h.*
  *Further, if $\pi^T G = h\pi^T$ for some h then $h = 1 + \gamma$.*

Proof: (a) Suppose $Ge = ge$ for some g. Then, from Theorem 2, $A \equiv I - (I - P)G = \alpha\pi^T$ so that $\alpha = Ae = [I - (I - P)G]e = e - g(I - P)e = e$, since $Pe = e$. i.e. $\alpha = e$ implying $G \in A\{1, 5a\}$.

Conversely, suppose $\alpha = e$, and let $Ge = g$. Then $e = Ae \equiv [I - (I - P)G]e = e - (I - P)g$ implying $(I - P)g = 0$, or that $g = Pg$. Thus, since $g$ is a right eigenvector of $P$ corresponding to the eigenvalue 1 for a finite irreducible Markov chain, $g = ge$ for some g. (For a proof see Theorem 6.1.5 of Hunter, (1983) for regular chains and Theorem 6.1.6 of Hunter, (1983) for periodic chains).

Further if $\boldsymbol{\alpha} = \boldsymbol{e}$, from (2.7), we must have $G\boldsymbol{e} = (1 + \gamma)\boldsymbol{e} = g\boldsymbol{e}$, establishing the observation that $g = 1 + \gamma$.

(b) Suppose $\boldsymbol{\pi}^T G = h\boldsymbol{\pi}^T$ for some $h$. Then, from Theorem 2, $B \equiv I - G(I - P) = \boldsymbol{e}\boldsymbol{\beta}^T$ so that $\boldsymbol{\beta}^T = \boldsymbol{\pi}^T B = \boldsymbol{\pi}^T[I - G(I - P)] = \boldsymbol{\pi}^T - h\boldsymbol{\pi}^T(I - P) = \boldsymbol{\pi}^T$, since $\boldsymbol{\pi}^T = \boldsymbol{\pi}^T P$. i.e. $\boldsymbol{\beta} = \boldsymbol{\pi}$ implying $G \in A\{1, 5b\}$.

Conversely, suppose $\boldsymbol{\beta} = \boldsymbol{\pi}$, and let $\boldsymbol{\pi}^T G = \boldsymbol{h}^T$. Then $\boldsymbol{\pi}^T = \boldsymbol{\pi}^T B = \boldsymbol{\pi}^T[I - G(I - P)] = \boldsymbol{\pi}^T - \boldsymbol{h}^T(I - P)$, implying $\boldsymbol{h}^T(I - P) = \boldsymbol{0}^T$, or that $\boldsymbol{h}^T = \boldsymbol{h}^T P$. Thus, since P is the transition matrix of a finite irreducible Markov chain $\boldsymbol{h}^T$ is a multiple of the stationary probability vector so that $\boldsymbol{h}^T = h\boldsymbol{\pi}^T$ for some $h$.

Further if $\boldsymbol{\beta} = \boldsymbol{\pi}$, from (2.7), we must have $\boldsymbol{\pi}^T G = (\gamma + 1)\boldsymbol{\pi}^T = h\boldsymbol{\pi}^T$, establishing the fact that $h = 1 + \gamma$. □

**Corollary 3.1**: *Let G have the characterisation* $G = G(\boldsymbol{\alpha}, \boldsymbol{\beta}, \gamma)$. *If* $G\boldsymbol{e} = g\boldsymbol{e}$ *for some g, and if* $\boldsymbol{\pi}^T G = h\boldsymbol{\pi}^T$ *for some h, then* $g = h = \gamma + 1$.
*Further this implies that* $G = G(\boldsymbol{e}, \boldsymbol{\pi}, \gamma)$ *and* $G \in A\{1, 5\}$.

Proof: Follows directly from Theorem 3 since if $G\boldsymbol{e} = g\boldsymbol{e}$ for some g, $G \in A\{1, 5a\}$ implying $\boldsymbol{\alpha} = \boldsymbol{e}$. Also since $\boldsymbol{\pi}^T G = h\boldsymbol{\pi}^T$ $G \in A\{1, 5b\}$ implying $\boldsymbol{\beta} = \boldsymbol{\pi}$. The value of $\gamma$ must be the common value of $g - 1$ and $h - 1$. □

Well known unique g-inverses are $A^{(1,2,3,4)}$, the "Moore-Penrose" g-inverse, (Moore, (1920)), (Penrose, (1955)); and $A^{(1,2,5)}$, the "group inverse" (Erdelyi, (1967)).

In Markov chain theory the following special multi-condition g-inverses of $I - P$ have been considered.

(a) Kemeny and Snell's *fundamental matrix* of finite irreducible Markov chains, $Z = [I - P + \Pi]^{-1}$ where $\Pi = \boldsymbol{e}\boldsymbol{\pi}^T$ was introduced by Kemeny and Snell, (1960). Z was shown in Hunter, (1969) to be a one-condition g-inverse of $I - P$ but it is also a (1, 5) g-inverse, (Hunter, (1988)).

(b) The *group inverse* of $I - P$ is the unique matrix $A^{(1,2,5)}$ satisfying conditions 1, 2 and 5 is $A^{\#} \equiv [I - P + \Pi]^{-1} - \Pi$. $A^{\#}$ was originally identified as the group inverse of $I - P$ in Meyer, (1975).

(c) The *Moore-Penrose* g-inverse of $I - P$ is the unique matrix $A^{(1,2,3,4)}$ satisfying conditions 1, 2, 3 and 4 and can be expressed (Hunter, (1988)) as $A^{(1,2,3,4)} = [I - P + \boldsymbol{\pi}\boldsymbol{e}^T]^{-1} - \dfrac{\boldsymbol{e}\boldsymbol{\pi}^T}{m\boldsymbol{\pi}^T\boldsymbol{\pi}}$.

An equivalent form, $A^{(1,2,3,4)} = [I - P + \alpha\boldsymbol{\pi}\boldsymbol{e}^T]^{-1} - \alpha\boldsymbol{e}\boldsymbol{\pi}^T$, where $\alpha = 1/\sqrt{m\boldsymbol{\pi}^T\boldsymbol{\pi}}$ was originally derived by Paige, Styan and Wachter, (1975). The equivalence of the two expressions was established in Hunter, (1988).

## 3. Mean first passage times.

The following result is well known, see for example, Kemeny and Snell, (1960).

**Theorem 4**: *The mean first passage times $m_{ij}$, ($1 \leq i \leq m$, $1 \leq j \leq m$), satisfy the following equation*

$$m_{ij} = 1 + \sum_{k \neq j} p_{ik} m_{kj}. \tag{3.1}$$

*In particular, the mean recurrence time of state j is given by*

$$m_{jj} = 1/\pi_j. \tag{3.2}$$

*The mean first passage time matrix $M = [m_{ij}]$ satisfies the matrix equation*

$$(I - P)M = E - PM_d, \tag{3.3}$$

*where $E = ee^T = [1]$, $M_d = [\delta_{ij} m_{jj}]$, a diagonal matrix with elements the diagonal elements of M. so that, $M_d = (\Pi_d)^{-1} \equiv D$, where $\Pi = e\pi^T$.* □

Hunter, (1982) showed that the solution of equations of the form of (3.3) can be implemented using g-inverses of $I - P$.

**Theorem 5**: *If G is any g-inverse of $I - P$, then the mean first passage time matrix, M, can be expressed as*

$$M = [G\Pi - E(G\Pi)_d + I - G + EG_d]D. \tag{3.4}$$

*If $H \equiv G(I - \Pi)$ then*

$$M = [I - H + EH_d]D. \tag{3.5}$$

Proof: Equation (3.4) was derived in Hunter, (1982). The special case given by equation (3.5) appears in Hunter, (2006). □

The advantage of the alternative expression (3.5) is that it leads to simpler elemental forms for the $m_{ij}$. If $A = [a_{ij}]$ is a matrix, define $a_{i\bullet} = \sum_{j=1}^{m} a_{ij}$ and $a_{\bullet j} = \sum_{i=1}^{m} a_{ij}$. Let $\delta_{ij} = 1$, when $i = j$ and 0 when $i \neq j$.

**Corollary 5.1**: *If $G = [g_{ij}]$ and $H = [h_{ij}]$, then*
(i) $\quad m_{ij} = ([g_{jj} - g_{ij} + \delta_{ij}]/\pi_j) + (g_{i\bullet} - g_{j\bullet}),$ *for all i, j.* (3.6)
(ii) $\quad m_{ij} = ([h_{jj} - h_{ij} + \delta_{ij}]/\pi_j)$ *for all i, j.* (3.7)
□

Note that for each form, when $i = j$, $m_{jj} = 1/\pi_j$, as prescribed by equation (3.2).

The expression given by (3.4) requires the awkward computation of $G\Pi - E(G\Pi)_d$, while (3.6) requires the computations of the row sums of $G$, $g_{i\bullet}$, with related similar computations for the alternative expressions. However, for $G \in A\{1, 5a\}$ these expressions need not be computed. In fact we have the following result:

**Theorem 6**: $G \in A\{1, 5a\}$ *if and only if*
$$M = [I - G + EG_d]D. \tag{3.8}$$

*Further if* $G = [g_{ij}] \in A\{1, 5a\}$ *then*
$$m_{ij} = ([g_{jj} - g_{ij} + \delta_{ij}]/\pi_j) \text{ for all } i,j$$
$$= \begin{cases} 1/\pi_j, & i = j, \\ (g_{jj} - g_{ij})/\pi_j & i \neq j. \end{cases} \tag{3.9}$$

Proof: First note that if equation (3.8) holds, then from (3.4), $G\Pi = E(G\Pi)_d$. Equating the $(i,j)$th elements we obtain $g_{i\bullet}\pi_j = g_{j\bullet}\pi_j$, so that for all $i, j$, $g_{i\bullet} = g_{j\bullet}$. This implies that $g_{i\bullet} = e_i Ge = g$, is a constant for all $i$, so that the sum of the elements of each row of $G$ is the same i.e. $Ge = ge$. Thus from Theorem 3(a), $G \in A\{1, 5a\}$.

Conversely, if $G \in A\{1, 5a\}$ then $Ge = ge$ for some $g$. In Hunter, (2007) it was shown that under any of the following three equivalent conditions: *(i) $Ge = ge$, $g$ a constant, (ii) $GE - E(G\Pi)_d D = 0$, (iii) $G\Pi - E(G\Pi)_d = 0$*; that (3.4) simplifies to (3.8).  □

Note that equation (3.5) of Theorem 5 is of the same form as equation (3.8) of Theorem 6 indicating that $H$ is a member of $A\{1, 5a\}$. We can in fact deduce more:

**Theorem 7**: *Let $G$ be any one-condition g-inverse of $I - P$ with characterisation $G = G(\alpha, \beta, \gamma)$. Let $H = G(I - \Pi)$. Then $H$ is a one-condition g-inverse of $I - P$, and $H$ has the characterization $H = G(e, \beta, -1)$ and hence is a member of the g-inverse class $A\{1, 2, 5a\}$.*

Proof: $H$ will be a 1-condition g-inverse of $I - P$ if $(I - P)H(I - P) = (I - P)$. From (2.5) of Theorem 2, $(I - P)G = I - A = I - \alpha\pi^T$. Substitution for $H = G(I - \Pi)$ and simplification leads to the required verification.
Let us suppose that $H$ has g-inverse classification parameters $\alpha_H, \beta_H, \gamma_H$. From Theorem 2 if $A_H = I - (I - P)H$ then $\alpha_H = A_H e$. Thus $\alpha_H = [I - (I - P)H]e = e - (I - P)G(I - \Pi)]e =$
$= e - (I - \alpha\pi^T)(I - e\pi^T)e$ which simplifies to yield $\alpha_H = e$. Also from Theorem 2, if $B_H = I - H(I - P)$ then $\beta_H^T = \pi^T B_H$. Thus $\beta_H^T = \pi^T[I - H(I - P)] = \pi^T - \pi^T G(I - \Pi)(I - P)$
i.e. $\beta_H^T = \pi^T - \pi^T G(I - P)$. From (2.5) and (2.6), $\beta_H^T = \pi^T B = \beta^T$.
Further, from (2.8), $\gamma + 1 = \beta^T G\alpha$, so that $\gamma_H + 1 = \beta_H^T H\alpha_H = \beta^T G(I - \Pi)e =$
$\beta^T G(I - e\pi^T)e = \beta^T Ge - \beta^T Ge\pi^T e = \beta^T Ge - \beta^T Ge = 0$ implying $\gamma_H = -1$.  □

Theorem 6 holds for any $G$ characterised as $G = G(e, \beta, \gamma)$. Well-known special cases for equation (3.8) are $G = Z = G(e, \pi, 0)$, Kemeny and Snell's fundamental matrix of the irreducible Markov chain, and $G = A^\# = G(e, \pi, -1)$, Meyer's group inverse of $I - P$.

Further, if $Z = [z_{ij}]$ and $A^\# = [a_{ij}^\#]$, equation (3.9) holds with either $g_{ij} = z_{ij}$ or $g_{ij} = a_{ij}^\#$.

Note that $H$ has some of the characteristics ($\alpha = e$ and $\gamma = -1$) of $A^\#$ but not all. One additional specification will in fact characterise $A^\#$: If $G$ is any one condition g-inverse of $I - P$ then $(I - \Pi)G(I - \Pi) = (I - \Pi)H = A^{\#}$. This result is well known, having appeared in Theorem 6.3 of Hunter, (1982).

We explore the simplification of the general equation (3.4) under the characterisation conditions of $G$, as described in Theorem 2.

**Theorem 8**: *Let $G$ have the characterisation $G = G(\alpha, \beta, \gamma)$ and define $Ge = g$ and $H = G(I - \Pi)$. Then*

$$M\Pi_d \alpha = \alpha + g + e[e^T H_d \alpha - (1+\gamma)], \tag{3.10}$$

$$\beta^T M\Pi_d = \beta^T + e^T H_d. \tag{3.11}$$

*In particular,* $ee^T H_d \alpha = (M\Pi_d - I)\alpha - g + e(1+\gamma) = e\beta^T(M\Pi_d - I)\alpha$. (3.12)

Proof: From equation (3.4),
$$M\Pi_d = G\Pi - E(G\Pi)_d + I - G + EG_d.$$

From (2.4) and (2.7), $\pi^T \alpha = 1$ and $G\alpha = (\gamma + 1)e$, so that

$$M\Pi_d \alpha = Ge\pi^T \alpha - ee^T(Ge\pi^T)_d \alpha + \alpha - G\alpha + ee^T G_d \alpha$$
$$= g - ee^T(g\pi^T)_d \alpha + \alpha - (\gamma + 1)e + ee^T G_d \alpha$$
$$= \alpha + g + e[e^T(G - G\Pi)_d \alpha - (1+\gamma)], \text{ leading to (3.10).}$$

Further, from (2.4) and (2.7), $\beta^T e = 1$ and $\beta^T G = (\gamma + 1)\pi^T$, so that

$$\beta^T M\Pi_d = \beta^T Ge\pi^T - \beta^T ee^T(Ge\pi^T)_d + \beta^T - \beta^T G + \beta^T ee^T G_d$$
$$= (\gamma+1)\pi^T - e^T(g\pi^T)_d + \beta^T - (\gamma+1)\pi^T + e^T G_d$$
$$= \beta^T + e^T(G - G\Pi)_d, \text{ leading to (3.11).}$$

Equation (3.12) follows directly from (3.10) and (3.11). □

Our quest would be accomplished if we could find a matrix expression for $G$ in terms of the parameters $\alpha$, $\beta$ and $\gamma$ of the g-inverse and the stationary probability vector $\pi$ and the matrix $M$. While Theorem 9 gives some relationships we have not been able to deduce such matrix expressions. Consequently we now explore, in the section to follow, element-wise expressions.

## 4. G-inverses in terms of stationary probabilities and mean first passage times.

The key results that we have at our disposal to find expressions for the elements of the g-inverses in terms of the properties of the underlying Markov chain are expressions (2.7) and Theorem 5 and its corollary.

**Theorem 9**: *Let $G = [g_{ij}] = G(\alpha, \beta, \gamma)$ be any g-inverse of $I - P$, then*

$$(i) \quad g_{jj} = \frac{1}{\alpha_j}\left(1 + \gamma - \sum_{k \neq j} \alpha_k g_{jk}\right). \tag{4.1}$$

(ii) $$g_{jj} = \frac{1}{\beta_j}\left((1+\gamma)\pi_j - \sum_{i \neq j}\beta_i g_{ij}\right). \quad (4.2)$$

Proof: (i) From the first expression of (2.7) that, for all g-inverses $G$ of $I - P$, $G\boldsymbol{\alpha} = (1 + \gamma)\boldsymbol{e}$, we have that, for all $j$, $\sum_{k=1}^{m}\alpha_k g_{jk} = 1+\gamma$, leading to equation (4.1).

(ii) From the second expression of (2.7) that, for all g-inverses $G$ of $I - P$, $\boldsymbol{\beta}^T G = (\gamma + 1)\boldsymbol{\pi}^T$, we have that, for all $j$, $\sum_{k=1}^{m}\beta_k g_{jk} = (1+\gamma)\pi_j$, leading to equation (4.2). □

Theorem 9 suggests that if we have expressions for the $g_{ij}$ elements when $i \neq j$ we can deduce an expressions for $g_{jj}$. Note that

$$\sum_{j=1}^{m}\sum_{k=1}^{m}\alpha_k g_{jk} = \sum_{k=1}^{m}\alpha_k \sum_{j=1}^{m}g_{jk} = \sum_{k=1}^{m}\alpha_k g_{\bullet k} = m(1+\gamma). \quad (4.3)$$

Equations (4.1), (4.2) and (4.3) do not, in general, yield any information about $g_{i\bullet} = \sum_{j=1}^{m}g_{ij}$ for specific $i$, which is required if we are to use the generality of equation (3.6).

However, from equation (4.1) we can deduce an expression for $g_{i\bullet}$ if $\alpha_k = \alpha$, a constant ( = 1, since $1 = \boldsymbol{\pi}^T\boldsymbol{\alpha} = \alpha\boldsymbol{\pi}^T\boldsymbol{e} = \alpha$) leading to $g_{i\bullet} = 1+\gamma$. Thus in this special case, with $\boldsymbol{\alpha} = \boldsymbol{e}$, i.e. $G \in A\{1, 5a\}$, equation (3.6) reduces to equation (3.9). Note that from equation (3.6), when $g_{i\bullet} = g_{j\bullet}$, $g_{ij} = g_{jj} + \pi_j m_{ij}$ for all $i \neq j$, so that if we can find a suitable expression for $g_{jj}$, using say equation (4.2), we can find a general expression for all $g_{ij}$. We defer consideration of this special case when $\boldsymbol{\alpha} = \boldsymbol{e}$ to later (See Corollary 11.1).

In order to use equation (3.6) to find expressions for the individual $g_{ij}$, for an arbitrary g-inverse it is clear that some knowledge of the row sums $g_{i\bullet}$ of $G$ is required. Note $g_{i\bullet} = \boldsymbol{e}_i^T G\boldsymbol{e}$. Consequently we typically require knowledge of $G\boldsymbol{e} = \boldsymbol{g}$. As we have seen, simplification, leading to equation (3.9), occurs only when $\boldsymbol{g} = g\boldsymbol{e}$, for some $g$. We need an alternative approach that yields expressions for the $g_{i\bullet}$.

**Theorem 10**: If $G = G(\boldsymbol{\alpha}, \boldsymbol{\beta}, \gamma)$, is any g-inverse of $I - P$ and $H \equiv G(I - \Pi)$ then

(i) $\quad \boldsymbol{\beta}^T H = \boldsymbol{0}^T,$ (4.4)

(ii) $\quad H\boldsymbol{\alpha} = \left[(\gamma + 1)I - G\right]\boldsymbol{e}.$ (4.5)

Proof: (i) Since, from (2.7), $\boldsymbol{\beta}^T G = (\gamma + 1)\boldsymbol{\pi}^T$ observe that
$$\boldsymbol{\beta}^T H = \boldsymbol{\beta}^T G - \boldsymbol{\beta}^T G\boldsymbol{e}\boldsymbol{\pi}^T = (\gamma + 1)\boldsymbol{\pi}^T - (\gamma + 1)\boldsymbol{\pi}^T\boldsymbol{e}\boldsymbol{\pi}^T = (\gamma + 1)\boldsymbol{\pi}^T - (\gamma + 1)\boldsymbol{\pi}^T = \boldsymbol{0}^T.$$

(ii) Since, from (2.7), $G\boldsymbol{\alpha} = (\gamma+1)\boldsymbol{e}$ we have that
$$H\boldsymbol{\alpha} = G(I - \Pi)\boldsymbol{\alpha} = G\boldsymbol{\alpha} - G\boldsymbol{e}\boldsymbol{\pi}^T\boldsymbol{\alpha} = (\gamma + 1)\boldsymbol{e} - G\boldsymbol{e} = \left[(\gamma + 1)I - G\right]\boldsymbol{e}. \quad □$$

We are now in the position to express $h_{ij}$ in terms of the $\{\beta_k\}$ parameters of the g-inverse.

**Corollary 10.1**: If $H = [h_{ij}]$ and $\delta_j \equiv \sum_{k \neq j}^{m}\beta_k m_{kj}$ then

$$h_{jj} = \pi_j \delta_j, \; (j = 1, ...,m) \; and \; h_{ij} = \pi_j(\delta_j - m_{ij}), \; i \neq j, \; (i, j = 1, ...,m). \quad (4.6)$$

Proof: From equation (4.4),

$$0 = \sum_{i=1}^{m}\beta_i h_{ij} = \beta_j h_{jj} + \sum_{i\neq j}\beta_i h_{ij} = \beta_j h_{jj} + \sum_{i\neq j}\beta_i(h_{jj} - m_{ij}\pi_j),$$

$$= \beta_j h_{jj} + h_{jj}\sum_{i\neq j}\beta_i - \sum_{i\neq j}\beta_i m_{ij}\pi_j = h_{jj}\left(\sum_{i=1}^{m}\beta_i\right) - \pi_j\sum_{i\neq j}\beta_i m_{ij}$$

$$= h_{jj} - \pi_j\sum_{i\neq j}\beta_i m_{ij}, \text{ since, by (2.4), } \sum_{i=1}^{m}\beta_i = \boldsymbol{\beta}^T\boldsymbol{e} = 1, \text{ leading to (4.6) when } i = j.$$

The result for $i \neq j$ follows from equation (3.7), since $h_{ij} = h_{jj} - m_{ij}\pi_j$. □

Note that we earlier showed that $H$ is a member of the $A\{1, 2, 5a\}$ class of g-inverses so that the elemental forms for the $h_{ij}$ should also hold for g-inverses in that class. We verify this later (Corollary 11.2).

We can now find an expression for the row sums of $G$, viz $g_{i\bullet}$ in terms of the $\{\alpha_k\}$ parameters of the g-inverse and the $\pi_j$ and the $m_{ij}$.

**Corollary 10.2**: If $\delta_j \equiv \sum_{k\neq j}\beta_k m_{kj}$ then

$$g_{i\bullet} = 1 + \gamma + \sum_{k\neq i}\pi_k\alpha_k m_{ik} - \sum_{k=1}^{m}\pi_k\alpha_k\delta_k, \quad (j = 1, ..., m). \tag{4.7}$$

Proof: Since $\boldsymbol{\alpha}^T = (\alpha_1, \alpha_2, ..., \alpha_m)$, from (4.5), $\sum_{k=1}^{m} h_{ik}\alpha_k = 1 + \gamma - \sum_{k=1}^{m} g_{ik} = 1 + \gamma - g_{i\bullet}$, implying that $g_{i\bullet} = 1 + \gamma - \sum_{k=1}^{m} h_{ik}\alpha_k = 1 + \gamma - h_{ii}\alpha_i - \sum_{k\neq i} h_{ik}\alpha_k$. Using (4.6) this leads to

$g_{i\bullet} = 1 + \gamma - \pi_i\alpha_i\delta_i - \sum_{k\neq i}\pi_k\alpha_k(\delta_k - m_{ik})$, from which (4.7) follows. □

If we can now express $g_{ij}$ in terms of $h_{ij}$ we can find an expression for the elements of any g-inverse of $I - P$. This leads to the following general theorem:

**Theorem 11**: Let $G = G(\boldsymbol{\alpha}, \boldsymbol{\beta}, \gamma)$ be any g-inverse of $I - P$ then the elements of $G = [g_{ij}]$ can be expressed in terms of the parameters $\{\alpha_j\}, \{\beta_j\}, \gamma$, the stationary probabilities $\{\pi_j\}$, and the mean first passage times $\{m_{ij}\}$, of the Markov chain, with $\delta_j \equiv \sum_{k\neq j}\beta_k m_{kj}$ as

$$g_{ij} = \begin{cases} \left(1 + \gamma + \delta_j - m_{ij} + \sum_{k\neq i}\pi_k\alpha_k m_{ik} - \sum_{k=1}^{m}\pi_k\alpha_k\delta_k\right)\pi_j, & i \neq j, \\ \left(1 + \gamma + \delta_j + \sum_{k\neq j}\pi_k\alpha_k m_{jk} - \sum_{k=1}^{m}\pi_k\alpha_k\delta_k\right)\pi_j, & i = j. \end{cases} \tag{4.8}$$

Proof: From the definition of $H$ as $G(I - \Pi)$, $h_{ij} = g_{ij} - g_{i\bullet}\pi_j$ implying that $g_{ij} = h_{ij} + g_{i\bullet}\pi_j$ for all $i, j$ which leads to (4.8) using (4.6) and (4.7) for the separate cases of $i = j$ and $i \neq j$. □

An interesting observation from equations (4.7) and (4.8) is that we can find some new interconnections:

For all one-condition g-inverses of $I - P$,

$$g_{ij} = \begin{cases} \pi_j(\delta_j + g_{i\bullet} - m_{ij}), & i \neq j, \\ \pi_j(\delta_j + g_{j\bullet}), & i = j, \end{cases} \qquad (4.9)$$

leading to an alternative expression for the $\delta_j = \dfrac{g_{jj}}{\pi_j} - g_{j\bullet}$, $j = 1, 2, \ldots, m$.

Theorem 11 is our key result giving new general expressions for the entire family of one condition g-inverses of $I - P$, in terms of the properties of the Markov chain and the classification parameters of the generalised inverse.

Before we consider some special cases we examine some simplification properties.

**Theorem 12**: *For all $G = G(\boldsymbol{\alpha}, \boldsymbol{\beta}, \boldsymbol{\gamma})$ g-inverses of $I - P$, with $\delta_j \equiv \sum_{k \neq j} \beta_k m_{kj}$ then*

$$\sum_{k=1}^{m} \pi_k \delta_k = K - 1, \qquad (4.10)$$

where $K = \sum_{j=1}^{m} \pi_j m_{ij}$, a constant (Kemeny's constant) for all $j = 1, \ldots, m$.

Proof: From the definition of $\delta_k$ and property (2.4) that $\sum_{i=1}^{m} \beta_i = 1$,

$$\sum_{k=1}^{m} \pi_k \delta_k = \sum_{k=1}^{m} \pi_k \sum_{i \neq k} \beta_i m_{ik} = \sum_{k=1}^{m} \pi_k \left( \sum_{i=1}^{m} \beta_i m_{ik} - \beta_k m_{kk} \right),$$

$$= \sum_{i=1}^{m} \beta_i \sum_{k=1}^{m} \pi_k m_{ik} - \sum_{k=1}^{m} \beta_k, \text{ (since } \pi_k m_{kk} = 1\text{),}$$

$$= (\sum_{i=1}^{m} \beta_i)(K - 1) = K - 1.$$

An alternative proof can be given by finding an expression from (4.9) for $g_{i\bullet} = \sum_{j=1}^{m} g_{ij}$ leading to $\sum_{j=1}^{m} \pi_j m_{ij} = \sum_{j=1}^{m} \pi_j \delta_j + 1$, from which (4.10) follows. □

(See Section 7 for a further discussion on Kemeny's constant.)

**Corollary 11.1**: *If $G$ has the characterisation $G \in A\{1, 5a\}$, i.e. $G = G(\boldsymbol{e}, \boldsymbol{\beta}, \gamma)$, with $\boldsymbol{\beta}^T = (\beta_1, \beta_2, \ldots, \beta_m)$, then the elements of $G = [g_{ij}]$ can be expressed in terms of the parameters $\{\beta_j\}$, $\gamma$, the stationary probabilities $\{\pi_j\}$, and the mean first passage times $\{m_{ij}\}$, of the Markov chain as*

$$g_{ij} = \begin{cases} \pi_j(1 + \gamma - m_{ij} + \delta_j), & i \neq j, \\ \pi_j(1 + \gamma + \delta_j), & i = j, \end{cases} \qquad (4.11)$$

where $\delta_j \equiv \sum_{k \neq j}^{m} \beta_k m_{kj}$.

Proof: When $\boldsymbol{\alpha} = \boldsymbol{e}$, from equation (4.8),

$$g_{ij} = \begin{cases} \left(1 + \gamma + \delta_j - m_{ij} + \sum_{k \neq i} \pi_k m_{ik} - \sum_{k=1}^{m} \pi_k \delta_k \right)\pi_j, & i \neq j, \\ \left(1 + \gamma + \delta_j + \sum_{k \neq j} \pi_k m_{jk} - \sum_{k=1}^{m} \pi_k \delta_k \right)\pi_j, & i = j, \end{cases}$$

which, using Theorem 12, yields equation (4.11).

(Note that, for $G \in A\{1, 5a\}$, equation (4.11) can be derived directly from equation (4.2) by substituting for $i \neq j$, $g_{ij} = g_{jj} - m_{ij}\pi_j$, (from equation (3.9)). Simplification leads directly to equation (4.11) for $g_{jj}$ and consequently for $g_{ij}$.)  □

We now consider a variety of special classes of generalized inverses, initially with the restriction that $\alpha = e$, i.e. $G \in A\{1, 5a\}$.

**Corollary 11.2**: *Let $G \in A\{1, 2, 5a\}$ and have the characterisation $G = G(e, \beta, -1)$. Then the elements of $G = [g_{ij}]$ can be expressed in terms of the $\{\beta_j\}$, the stationary probabilities $\{\pi_j\}$, and the mean first passage times $\{m_{ij}\}$, of the Markov chain as*

$$g_{ij} = \begin{cases} \pi_j(\delta_j - m_{ij}), & i \neq j, \\ \pi_j \delta_j, & i = j, \end{cases} \quad (4.12)$$

*where $\delta_j = \sum_{k \neq j}^{m} \beta_k m_{kj}$.*

Proof: From (2.9), $G \in A\{1, 2\} \Leftrightarrow \gamma = -1$, and thus (4.12) follows from (4.11). Note also that $G$ has the property that $Ge = 0$, implying that $g_{i\cdot} = 0$ for all $i$. Note also that $\pi^T G = 0^T$.  □

Note that the conclusion of Corollary 11.2 is consistent with the observation made following Corollary 10.1 in respect to the elements of $H = [h_{ij}]$.

**Corollary 11.3**: *Let $G \in A\{1, 4, 5a\}$ and have the characterisation $G = G(e, e/m, \gamma)$. Then the elements of $G = [g_{ij}]$ can be expressed in terms of $\gamma$, the stationary probabilities $\{\pi_j\}$, and the mean first passage times $\{m_{ij}\}$, of the Markov chain with $\eta_j \equiv (\sum_{k \neq j}^{m} m_{kj})/m$, as*

$$g_{ij} = \begin{cases} \pi_j(1+\gamma+\eta_j - m_{ij}), & i \neq j, \\ \pi_j(1+\gamma+\eta_j), & i = j. \end{cases} \quad (4.13)$$

Proof: From (2.11), $G \in A\{1, 4\} \Leftrightarrow \beta = e/m$, and thus (4.13) follows from (4.11) with $\delta_j$ replaced by $\eta_j$. Note also that $e^T G = m(\gamma+1)\pi^T$.  □

**Corollary 11.4**: *Let $G \in A\{1, 2, 4, 5a\}$ and have the characterisation $G = G(e, e/m, -1)$. Then the elements of $G = [g_{ij}]$ can be expressed in terms of the stationary probabilities $\{\pi_j\}$, and the mean first passage times $\{m_{ij}\}$, of the Markov chain, with $\eta_j \equiv (\sum_{k \neq j}^{m} m_{kj})/m$, as*

$$g_{ij} = \begin{cases} \pi_j(\eta_j - m_{ij}), & i \neq j, \\ \pi_j \eta_j, & i = j. \end{cases} \quad (4.14)$$

Proof: Equation (4.14) follows from (4.13) with $\gamma = -1$ or (4.12) with $\beta_i = 1/m$ and $\delta_j$ replaced by $\eta_j$. Note also that $Ge = 0$ and $e^T G = 0^T$, implying that $g_{i\cdot} = 0$ and $g_{\cdot j} = 0$ for all $i$ and $j$.  □

Let us define, for each $j = 1, 2, \ldots, m$, $\tau_j \equiv \sum_{k=1}^{m} \pi_k m_{kj} = \sum_{k \neq j} \pi_k m_{kj} + 1$, (since $\pi_j m_{jj} = 1$.)

Thus if $G \in A\{1, 5b\}$ with $\boldsymbol{\beta} = \boldsymbol{\pi}$, then, in Theorem 11, $\delta_j = \tau_j - 1$, leading to appropriate modifications in equation (4.8). However, there is no significant simplification of the expressions for the $g_{ij}$.

**Corollary 11.5**: *Let $G \in A\{1, 5\}$ so that $G$ has the characterisation $G = G(\boldsymbol{e}, \boldsymbol{\pi}, \gamma)$. The elements of $G = [g_{ij}]$ can be expressed in terms of $\gamma$, the stationary probabilities $\{\pi_j\}$, and the mean first passage times $\{m_{ij}\}$, of the Markov chain as*

$$g_{ij} = \begin{cases} \pi_j(\tau_j + \gamma - m_{ij}), & i \neq j, \\ \pi_j(\tau_j + \gamma), & i = j. \end{cases} \quad (4.15)$$

Proof: From (2.12), $G \in A\{1, 5\} \Leftrightarrow \boldsymbol{\alpha} = \boldsymbol{e}, \boldsymbol{\beta} = \boldsymbol{\pi}$ and thus (4.15) follows from (4.11) bearing in mind that since $\boldsymbol{\alpha} = \boldsymbol{e}$, $G\boldsymbol{\alpha} = G\boldsymbol{e} = (\gamma + 1)\boldsymbol{e}$. Further $\boldsymbol{\pi}^T G = (1 + \gamma)\boldsymbol{\pi}^T$. □

**Corollary 11.6**: Let $Z = \left[z_{ij}\right] = [I - P + \boldsymbol{e}\boldsymbol{\pi}^T]^{-1}$ be Kemeny and Snell's fundamental matrix, then

$$z_{ij} = \begin{cases} \pi_j(\tau_j - m_{ij}), & i \neq j, \\ \pi_j \tau_j, & i = j. \end{cases} \quad (4.16)$$

Proof: $Z$ is an $A\{1, 5\}$ generalized inverse of $I - P$ with $\gamma = 0$ and (4.16) follows from (4.15). Note that $Z\boldsymbol{e} = \boldsymbol{e}$ and $\boldsymbol{\pi}^T Z = \boldsymbol{\pi}^T$. □

**Corollary 11.7**: Let $A^{\#} = \left[a_{ij}^{\#}\right] = [I - P + \boldsymbol{e}\boldsymbol{\pi}^T]^{-1} - \boldsymbol{e}\boldsymbol{\pi}^T$ be Meyer's group inverse of $I - P$, then

$$a_{ij}^{\#} = \begin{cases} \pi_j(\tau_j - 1 - m_{ij}), & i \neq j, \\ \pi_j(\tau_j - 1), & i = j. \end{cases} \quad (4.17)$$

Proof: Observe that $A^{\#}$ is the unique member of $A\{1, 2, 5\}$ and has the characterisation $G(\boldsymbol{e}, \boldsymbol{\pi}, -1)$ and (4.17) follows from (4.15) with $\gamma = -1$. Note also that $A^{\#}\boldsymbol{e} = \boldsymbol{0}$ and $\boldsymbol{\pi}^T A^{\#} = \boldsymbol{0}^T$. □

An expression equivalent to equation (4.17) also appears in Ben-Ari and Neumann, (2012), for discrete Markov chains where they attribute the result to Equation (3.6) of Meyer (1975) (which is actually equation (3.9) above with $g_{ij} = a_{ij}^{\#}$). While their result is correct, their conclusion requires knowledge of additional properties of $A^{\#}$. Equation (4.17) first requires deducing the elemental expression for $a_{jj}^{\#}$, which effectively follows from the facts that $A^{\#}\boldsymbol{e} = \boldsymbol{0}$ and $\boldsymbol{\pi}^T A^{\#} = \boldsymbol{0}^T$ (implying, effectively in our context, that $\boldsymbol{\alpha} = \boldsymbol{e}$, $\boldsymbol{\beta} = \boldsymbol{\pi}$ and $\gamma = -1$.)

The other main class of g-inverses that we have not considered involves the condition 3. Note that $G \in A\{1, 3\} \Leftrightarrow \boldsymbol{\alpha} = \boldsymbol{\pi}/\boldsymbol{\pi}^T \boldsymbol{\pi}$.

**Corollary 11.8**: *Let $G \in A\{1, 3\}$ so that G has the characterisation $G(\pi / \pi^T \pi, \beta, \gamma)$. Then the elements of $G = [g_{ij}]$ can be expressed in terms of the parameters $\{\alpha_j\}, \{\beta_j\}, \gamma$, the stationary probabilities $\{\pi_j\}$, and the mean first passage times $\{m_{ij}\}$, of the Markov chain, with $\delta_j \equiv \sum_{k \neq j}^{m} \beta_k m_{kj}$ as*

$$g_{ij} = \begin{cases} \left(1 + \gamma + \delta_j - m_{ij} + (1/\sum \pi_k^2)\sum_{k \neq i} \pi_k^2 m_{ik} - (1/\sum \pi_k^2)\sum_{k=1}^{m} \pi_k^2 \delta_k \right)\pi_j, & i \neq j, \\ \left(1 + \gamma + \delta_j + (1/\sum \pi_k^2)\sum_{k \neq j} \pi_k^2 m_{jk} - (1/\sum \pi_k^2)\sum_{k=1}^{m} \pi_k^2 \delta_k \right)\pi_j, & i = j. \end{cases} \quad (4.18)$$

□

This leads, with the additional restrictions of $\gamma = -1$, $\beta = e/m$ to a general expression for the elements of the Moore-Penrose generalised inverse of $I - P$.

**Corollary 11.9**: *Let $G \in A\{1, 2, 3, 4\}$, the Moore-Penrose g-inverse of $I - P$, so that G has the characterisation $G = G(\pi / \pi^T \pi, e/m, -1)$. The elements of $G = [g_{ij}]$ can be expressed in terms of the parameters $\{\alpha_j\}, \{\beta_j\}, \gamma$, the stationary probabilities $\{\pi_j\}$, and the mean first passage times $\{m_{ij}\}$, of the Markov chain, with $\eta_j \equiv (\sum_{k \neq j}^{m} m_{kj})/m$, as*

$$g_{ij} = \begin{cases} \left(\eta_j - m_{ij} + (1/\sum \pi_k^2)\sum_{k \neq i} \pi_k^2 m_{ik} - (1/\sum \pi_k^2)\sum_{k=1}^{m} \pi_k^2 \eta_k \right)\pi_j, & i \neq j, \\ \left(\eta_j + (1/\sum \pi_k^2)\sum_{k \neq j} \pi_k^2 m_{jk} - (1/\sum \pi_k^2)\sum_{k=1}^{m} \pi_k^2 \eta_k \right)\pi_j, & i = j. \end{cases} \quad (4.19)$$

□

As we have seen in the corollaries to Theorem 11, major simplification to expressions for *M* takes place only when $Ge = ge$ for some $g$, and this implies the necessity of the property, from Theorem 3, that $\alpha = e$. In the Moore-Penrose case, this is not possible unless $\pi^T = e^T/m$, in which case the Markov chain has a doubly stochastic transition matrix. See Hunter, (2010).

The results above substantiate the claims that the group inverse of $I - P$ contains a richness and simplicity of information about the stationary probabilities and the mean first passage times of the associated irreducible Markov chain.

## 5. Higher moments of the passage times

Let $m_{ij}^{(2)}$ be the *k-th* moment of the first passage time from state *i* to state *j* in the finite irreducible Markov chain with transition matrix *P* and let $M^{(2)} = [m_{ij}^{(2)}]$.

In Hunter, (2008) some new expressions, using g-inverses of $I - P$ were derived for $M^{(2)}$. Some of these expressions involving $M_d^{(2)} = [\delta_{ij} m_{jj}^{(2)}]$ are summarised below:

**Theorem 13:**
(*a*) $M^{(2)}$ *satisfies the matrix equation*
$$(I - P)M^{(2)} = E + 2P(M - M_d) - PM_d^{(2)}, \qquad (5.1)$$
*where*
$$M_d^{(2)} = 2D(\Pi M)_d - D, \qquad (5.2)$$
*with* $D = M_d = (\Pi_d)^{-1}$,

*implying* 
$$m_{jj}^{(2)} + m_{jj} = 2m_{jj}\sum_{i=1}^{m}\pi_i m_{ij}. \qquad (5.3)$$

(*b*) *If G is any g-inverse of I – P,* $\quad M_d^{(2)} = D + 2D\{(I - \Pi)G(I - \Pi)\}_d D.$ (5.4)

(*c*) *If* $G \in A\{1, 5a\}$ $\quad M_d^{(2)} = D + 2DG_d D - 2D(\Pi G)_d D.$ (5.5)

(*d*) *If* $G \in A\{1, 5b\}$ $\quad M_d^{(2)} = D + 2DG_d D - 2D(G\Pi)_d D.$ (5.6)

(*e*) *If* $G \in A\{1, 5\}$ $\quad M_d^{(2)} = 2DG_d D - (1 + 2\gamma)D.$ (5.7)

*In particular,* $\quad M_d^{(2)} = D + 2DA_d^\# D = 2DZ_d D - D.$ (5.8)

Proof:
(*a*) Equation (5.1) appears in Theorem 2.4 of Hunter, (2008), originally appearing in Theorem 4.5.1 of Kemeny and Snell, (1960). Equation (5.2) is given in Corollary 2.4.1, Hunter, (2008). Equation (5.3) appears in Corollary 2.5.5, Hunter, (2008).
(*b*) Equation (5.4) appears in Corollary 2.4.2, Hunter, (2008).
(*c*) Equation (5.5) appears in Corollary 2.4.3, Hunter, (2008).
(*d*) Equation (5.6) follows directly from equation (5.4) under the requisite simplification.
(*e*) Equation (5.7) follows directly from equation (5.6), since from Theorem 2,
$G\Pi = Ge\pi^T = (\gamma + 1)e\pi^T = (\gamma + 1)\Pi$ with $(G\Pi)_d = (\gamma + 1)\Pi_d = (\gamma + 1)D^{-1}$.
Equation (5.8) follows from (5.7) by noting that for $G = Z$, $\gamma = 0$ while for $G = A^\#$, $\gamma = -1$.
□

Note that if *G* is any g-inverse of *I – P*, the expression $(I - \Pi)G(I - \Pi)$ is invariant and is in fact the group inverse $A^\#$ (see Theorem 6.3, Hunter, (1982)). This observation is also reflected in Theorem 13 with equation (5.8) following directly from equation (5.4).

A consequence of Theorem 13 is that, from equation (5.3):
$$\tau_j = \sum_{i=1}^{m}\pi_i m_{ij} = \frac{m_{jj}^{(2)} + m_{jj}}{2m_{jj}} = \frac{1}{2} + \frac{m_{jj}^{(2)}}{2m_{jj}} = \frac{1 + \pi_j m_{jj}^{(2)}}{2}. \qquad (5.9)$$

In Corollary 2.5.6 and 2.5.7 of Hunter, (2008) more general expressions for $\tau_j$ are given in terms of g-inverses of *I – P* and summarised below:

**Theorem 14:**
(*a*) *If* $\tau^T = (\tau_1, \tau_2, ..., \tau_m)$ *then* $\tau^T = \pi^T M = e^T (\Pi M)_d$ *and* $\tau = (\Pi M)_d e.$ (5.10)

(*b*) *If G is any g-inverse of I – P, then* $\tau = e + (\pi^T Ge)e - (\Pi G)_d De - (GE)_d e + G_d De,$ (5.11)

*implying* $\tau_j = 1 + \sum_{i=1}^{m}\pi_i g_{i\bullet} - g_{j\bullet} + \left(g_{jj} - \sum_{i=1}^{m}\pi_i g_{ij}\right)\bigg/\pi_j.$ (5.12)

(*c*) *If* $G \in A\{1, 5a\}$ *then* $\tau = e - (\Pi G)_d De + G_d De,$ (5.13)

*implying* $\tau_j = 1 + \left(g_{jj} - \sum_{i=1}^{m}\pi_i g_{ij}\right)\Big/\pi_j.$ (5.14)

*(d) If* $G \in A\{1, 5b\}$ *then* $\boldsymbol{\tau} = \boldsymbol{e} - (GE)_d \boldsymbol{e} + G_d D\boldsymbol{e},$ (5.15)

*implying* $\tau_j = 1 - g_{j\bullet} + g_{jj}/\pi_j.$ (5.16)

*(e) If* $G \in A\{1, 5\}$ *then* $\boldsymbol{\tau} = G_d D\boldsymbol{e} - \gamma \boldsymbol{e},$ (5.17)

*implying, in particular,* $\boldsymbol{\tau} = Z_d D\boldsymbol{e} = \boldsymbol{e} + A_d^{\#} D\boldsymbol{e},$ *and* $\tau_j = (z_{jj}/\pi_j) = 1 + (a_{jj}^{\#}/\pi_j).$ (5.18)

Proof: Equations (5.10), (5.11), (5.13) appear in Corollary 2.5.6 of Hunter, (2008) while equations (5.12), (5.14) and (5.18) appear in Corollary 2.5.7 of Hunter, (2008). Equations (5.15), (5.16) and (5.17) follow from (5.11) and (5.12). □

The transpose variant of Theorem 14 (*e*), using *Z*, was initially derived by Kemeny and Snell, (1960) in their Theorem 4.49.

An immediate observation from Theorem 14(*e*) is that $z_{jj} = \pi_j \tau_j$ and $a_{jj}^{\#} = \pi_j(\tau_j - 1)$, as observed in Corollaries 11.6 and 11.7.

The crucial observation is that we have a variety of expressions for $\tau_j = \sum_{i=1}^{m}\pi_i m_{ij}.$

## 6. G-inverses in terms of stationary probabilities and the second moments of passage times.

Ben-Ari and Neumann, (2012), extend the property of equation (4.17) to show that for the group inverse $A^{\#} = \left[a_{ij}^{\#}\right]$, if $m_{ij}^{(2)} = E[T_{ij}^2 | X_0 = i]$ is the second moment of the passage time random variable $T_{ij}$, that

$$a_{ij}^{\#} = \begin{cases} \dfrac{1}{2}\left(\dfrac{m_{jj}^{(2)}}{(m_{jj})^2} - \dfrac{1}{m_{jj}}\right) - \pi_j m_{ij}, & i \neq j, \\ \dfrac{1}{2}\left(\dfrac{m_{jj}^{(2)}}{(m_{jj})^2} - \dfrac{1}{m_{jj}}\right), & i = j. \end{cases}$$ (6.1)

The proof that they give is based upon arguments involving analytic continuation, Laurent expansions and Taylor series expansions of generating functions. However a much simpler proof follows from Corollary 11.7 and equations (5.9). We do not need the full generality of the group inverse to get expressions for the elements of the generalized inverses involving the properties of the second moments of the first return time variable $T_{jj}$. The following result follows from Corollary 11.5 utilising equation (5.9).

**Theorem 15**: *Let* $G \in A\{1, 5\}$ *so that G has the characterisation* $G = G(\boldsymbol{e}, \boldsymbol{\pi}, \gamma)$. *The elements of* $G = [g_{ij}]$ *can be expressed in terms of* $\gamma$, *the stationary probabilities* $\{\pi_j\}$, *the mean first* passage times $\{m_{ij}\}$, and the second moment of the recurrence times $m_{jj}^{(2)}$, of the Markov chain as

$$g_{ij} = \begin{cases} \pi_j\left(\gamma + \dfrac{\pi_j m_{jj}^{(2)}+1}{2} - m_{ij}\right), & i \neq j, \\ \pi_j\left(\gamma + \dfrac{\pi_j m_{jj}^{(2)}+1}{2}\right), & i = j. \end{cases} \qquad (6.2)$$

□

Two special cases follow from Theorem 15, one with $\gamma = 0$ when $G = Z$, while the other with $\gamma = -1$ when $G = A^\#$ (as given by equation (6.1) above.)

## 7. Kemeny's constant.

Kemeny's constant is defined as $K = \sum_{j=1}^{m}\pi_j m_{ij} = \sum_{j \neq i}\pi_j m_{ij} + 1$ since $m_{ii} = 1/\pi_i$. The interesting observation is that this sum is in fact a constant, independent of $i$. This is in contrast to $\tau_j = \sum_{i=1}^{m}\pi_i m_{ij}$ which, as was shown in Section 5, varies with $j = 1, 2, \ldots, m$.

**Theorem 16**: *If $G = [g_{ij}]$ is any generalised inverse of $I - P$,*
$$K = 1 + \sum_{j=1}^{m}(g_{jj} - g_{j\bullet}\pi_j). \qquad (7.1)$$
Proof: From equations (3.6) observe that
$$K = 1 + \sum_{j \neq i}m_{ij}\pi_j = 1 + \sum_{j \neq i}(g_{jj} - g_{ij}) + \sum_{j \neq i}(g_{i\bullet} - g_{j\bullet})\pi_j,$$
$$= 1 + \sum_{j=1}^{m}(g_{jj} - g_{ij}) + \sum_{j=1}^{m}(g_{i\bullet} - g_{j\bullet})\pi_j = 1 + \sum_{j=1}^{m}g_{jj} - g_{i\bullet} + g_{i\bullet} - \sum_{j=1}^{m}g_{j\bullet}\pi_j,$$
leading to equation (7.1). □

Equation (7.1) was originally given in Hunter, (2006), (using a different notation) but not identified as Kemeny's constant until it was highlighted by Kirkland (2010).

**Theorem 17**: *If $G \in A\{1, 5a\}$, i.e. $G = G(e, \boldsymbol{\beta}, \gamma)$, then Kemeny's constant is given by*

$$K = tr(G) - \gamma. \qquad (7.2)$$
*In particular, $K = tr(Z) = tr(A^\#) + 1$.*

Proof: $G \in A\{1, 5a\}$ implies that $Ge = (1 + \gamma)e$ and thus $\sum_{j=1}^{m}g_{ij} = g_{i\bullet} = g = 1 + \gamma$.

Since $Ge = ge$, we can use equation (3.6). i.e. for $i \neq j$, $g_{ij} = g_{jj} - m_{ij}\pi_j$ implying

$$g = g_{i\bullet} = \sum_{j=1}^{m}g_{ij} = g_{ii} + \sum_{j \neq i}g_{ij} = g_{ii} + \sum_{j \neq i}(g_{jj} - m_{ij}\pi_j)$$
$$= \sum_{j=1}^{m}g_{jj} - \sum_{j \neq i}m_{ij}\pi_j = tr(G) - K + 1, \text{ leading to (7.2).} \qquad □$$

# 8. Applications to perturbed Markov chains

The ability to obtain specific expressions for the elements of particular g-inverses of $I - P$ leads us to explore bounds related to changes to the stationary distributions following perturbations of the transition matrix of the Markov chain. These bounds can, in special cases, be expressed in terms elements of the g-inverses and consequently in terms of the basic properties of the original Markov chain. viz. its stationary probabilities and its mean first passage times.

**Theorem 18:** *Let P be the transition matrix of a finite irreducible Markov chain. Let $\bar{P} = [\bar{p}_{ij}]$ be the transition matrix of the perturbed Markov chain so that $\bar{P} = P + \mathbf{E}$, where $\mathbf{E} = [\varepsilon_{ij}]$ is the matrix of perturbations (with $\sum_{j=1}^{m} \varepsilon_{ij} = 1$). We assume that $\bar{P}$ is irreducible. Let $\pi^T = (\pi_1, \pi_2, ..., \pi_m)$ and $\bar{\pi}^T = (\bar{\pi}_1, \bar{\pi}_2, ..., \bar{\pi}_m)$ be the stationary probability vectors of the Markov chains with transition matrices P and $\bar{P}$, respectively. Let $\Pi = e\pi^T$. Let G be any g-inverse of $I - P$ and let $H = G(I - \Pi)$. Then*

$$\bar{\pi}^T - \pi^T = \bar{\pi}^T E G(I - \Pi) = \bar{\pi}^T E H. \tag{8.1}$$

Proof: First observe that, since $\pi^T(I - P) = \mathbf{0}^T$ and $\bar{\pi}^T(I - \bar{P}) = \bar{\pi}^T(I - P - E) = \mathbf{0}^T$,

$$(\bar{\pi}^T - \pi^T)(I - P) = \bar{\pi}^T E. \tag{8.2}$$

Note, from equation (2.5) that $(I - P)G = I - \alpha\pi^T$, so that $(I - P)G(I - \Pi) = (I - \Pi)$. Further, note that $(\bar{\pi}^T - \pi^T)\Pi = (\bar{\pi}^T - \pi^T)e\pi^T = (\bar{\pi}^T e - \pi^T e)\pi^T = \mathbf{0}^T$. Thus post-multiplication of (8.2) by $H = G(I - \Pi)$ leads to equation (8.1). □

Theorem 18 appears as Theorem 2.1 in Hunter, (2005) and forms the basis for numerous expressions for the difference $\bar{\pi}^T - \pi^T$. In particular we can summarise the main matrix expressions as follows:

**Theorem 19**: *Under the conditions of Theorem 18, let G is any g-inverse of $I - P$ Let $M = [m_{ij}]$ be the mean first passage time matrix of the Markov chain associated with P. Let $\mathbf{E}$ be any general perturbation. Then*
*(i) If $G \in A\{1, 5a\}$, i.e. $G = G(e, \beta, \gamma)$, then $\bar{\pi}^T - \pi^T = \bar{\pi}^T E G$.* (8.3)

*(ii)* $\quad\quad\quad\quad\quad\quad\quad\quad\quad\quad\quad\quad \bar{\pi}^T - \pi^T = -\bar{\pi}^T E (M - M_d)(M_d)^{-1}.$ (8.4)

*(iii) If $N = (M - M_d)(M_d)^{-1}$,* $\quad\quad\quad \bar{\pi}^T - \pi^T = -\bar{\pi}^T E N.$ (8.5)

Proof:
(i) $G \in A\{1, 5a\}$, implies that $Ge = (1 + \gamma)e$ so that $\mathbf{E} H = \mathbf{E} G - \mathbf{E} Ge\pi^T = \mathbf{E} G - (1 + \gamma)\mathbf{E} e\pi^T = \mathbf{E} G$, since $\mathbf{E} e = \mathbf{0}$ and equation (8.3) follows from (8.1).
(ii) From (3.5) observe that $(M - M_d)(M_d)^{-1} = EH_d - H$ so that $H = EH_d - (M - M_d)(M_d)^{-1}$. Now $\mathbf{E} G(I - \Pi) = \mathbf{E} EH_d - \mathbf{E}(M - M_d(M_d)^{-1} = -\mathbf{E}(M - M_d)(M_d)^{-1}$ since $\mathbf{E} EH_d = \mathbf{E} ee^T H_d = 0$, and (8.4) follows from (8.1). (This was also given as Theorem 2.3, Hunter, (2005).)
(iii) Equation (8.5) follows directly from (8.4) since $N = (M - M_d)(M_d)^{-1}$. Alternatively use (8.1) and the fact that $N = EH_d - H$. (This result was also given as Theorem 2.4, Hunter, (2005)). □

Note that only when $G = G(e, \beta, \gamma),$ can we have a form of $\bar{\pi}^T - \pi^T$ involving just the generalized inverse $G$, as in (8.3). Special cases of (8.3) include

(a)   $G = [I - P + eu^T]^{-1} + ef^T + g\pi^T$ with $\pi^T t \neq 0$, $u^T e \neq 0$, $f^T$ and $g$ arbitrary vectors,
       leading to $\bar{\pi}^T - \pi^T = \bar{\pi}^T E [I - P + eu^T]^{-1}$.
(b)   $G = [I - P + eu^T]^{-1} + ef^T$ with $u^T e \neq 0$, and $f^T$ an arbitrary vector.
(c)   $G = Z$, the fundamental matrix of $I - P$.
(d)   $G = A^{\#}$, the group inverse of $I - P$.

The results (*a*) – (*d*) all appear in Hunter (2005). Result (*c*) was initially derived by Schweitzer, (1968). Result (*d*) is due to Meyer (1980). A special case of results (*a*) and (*b*), (with $f^T = 0^T$ and $g = 0$), appears in Seneta, (1988) while result (*b*) appears in Seneta, (1991).

We complete our discussion on applications of generalized inverses to investigate some element-wise expressions for the differences between the stationary probabilities $\pi_j$ and $\bar{\pi}_j$, utilising Theorem 18 and the element wise expressions that we have developed in section 7 of this paper.

**Theorem 20:** *Under the conditions of Theorem 18, let $G = [g_{ij}]$ be any generalized inverse of $I - P$, let $H = [h_{ij}] = G(I - \Pi)$, $M = [m_{ij}]$ $N = [n_{ij}]$, Then, for all $j = 1, 2, \ldots, m$,*

(i) $\qquad\qquad\qquad \bar{\pi}_j - \pi_j = \sum_{i=1}^{m}\sum_{k=1}^{m} \bar{\pi}_i \varepsilon_{ik} h_{kj}.$ (8.6)

(ii) $\qquad\qquad\qquad \bar{\pi}_j - \pi_j = \sum_{i=1}^{m}\sum_{k=1}^{m} \bar{\pi}_i \varepsilon_{ik} g_{kj} - \pi_j \sum_{i=1}^{m}\sum_{k=1}^{m} \bar{\pi}_i \varepsilon_{ik} g_{k\bullet}.$ (8.7)

(iii) If $G \in A\{1, 5a\}$, i.e. $G = G(e, \beta, \gamma)$,

$$\bar{\pi}_j - \pi_j = \sum_{i=1}^{m}\sum_{k=1}^{m} \bar{\pi}_i \varepsilon_{ik} g_{kj}. \qquad (8.8)$$

(iv) For all perturbations, $\pi_j - \bar{\pi}_j = \pi_j \sum_{i=1}^{m}\sum_{k \neq j} \bar{\pi}_i \varepsilon_{ik} m_{kj}.$ (8.9)

(v) For all perturbations, $\pi_j - \bar{\pi}_j = \sum_{l \neq j}\sum_{k=1}^{m} \bar{\pi}_k \varepsilon_{kl} n_{lj}.$ (8.10)

Proof:
(*i*) Equation (8.6) follows directly from (8.1).
(*ii*) Equation (8.7) follows from (8.7) since $h_{ij} = g_{ij} - g_{i\bullet}\pi_j$.
(*iii*) Equation (8.8) follows directly from (8.3) or from (8.7) by noting that the term involving $g_{k\bullet}$ disappears since in this situation $g_{k\bullet}$ is a constant and $\sum_{k=1}^{m} \varepsilon_{ik} = 0$.
(*iv*) Substitution in equation (8.1), using the expressions for $h_{kj}$ from equation (4.6), with $\delta_j \equiv \sum_{k \neq j}^{m} \beta_k m_{kj}$, yields

$$\bar{\pi}_j - \pi_j = \sum_{i=1}^{m}\sum_{k=1}^{m} \bar{\pi}_i \varepsilon_{ik} h_{kj} = \sum_{i=1}^{m} \bar{\pi}_i \left( \varepsilon_{ij} h_{jj} + \sum_{k \neq j}^{m} \varepsilon_{ik} h_{kj} \right)$$

$$= \sum_{i=1}^{m} \bar{\pi}_i \left( \varepsilon_{ij} \pi_j \delta_j + \sum_{k \neq j}^{m} \varepsilon_{ik} \pi_j (\delta_j - m_{kj}) \right)$$

$$= \pi_j \delta_j \sum_{i=1}^{m} \bar{\pi}_i \varepsilon_{ij} + \sum_{i=1}^{m} \bar{\pi}_i \left( \sum_{k \neq j}^{m} \varepsilon_{ik} \pi_j \delta_j \right) - \sum_{i=1}^{m} \bar{\pi}_i \left( \sum_{k \neq j}^{m} \varepsilon_{ik} \pi_j m_{kj} \right)$$

$$= \pi_j \delta_j \sum_{i=1}^{m} \bar{\pi}_i \varepsilon_{ij} + \pi_j \delta_j \sum_{i=1}^{m} \bar{\pi}_i \left( \sum_{k \neq j}^{m} \varepsilon_{ik} \right) - \pi_j \sum_{i=1}^{m} \bar{\pi}_i \left( \sum_{k \neq j}^{m} \varepsilon_{ik} m_{kj} \right)$$

$$= \pi_j \delta_j \sum_{i=1}^{m} \bar{\pi}_i \left( \sum_{k=1}^{m} \varepsilon_{ik} \right) - \pi_j \sum_{i=1}^{m} \bar{\pi}_i \left( \sum_{k \neq j}^{m} \varepsilon_{ik} m_{kj} \right)$$

$$= -\pi_j \sum_{i=1}^{m} \bar{\pi}_i \left( \sum_{k \neq j}^{m} \varepsilon_{ik} m_{kj} \right), \text{ since } \sum_{k=1}^{m} \varepsilon_{ik} = 0.$$

Note that an alternative derivation for (8.9) follows directly from an elemental expression of equations (8.4).

(*v*) Follows from (8.5) or using the property that $n_{ij} = (1 - \delta_{ij})\pi_j m_{ij}$. This was also given in Theorem 3.1 of Hunter, (2005). □

Note that substitution from Corollary 9.2 for $g_k$. and Theorem 11 for $g_{ij}$ in equation (8.7) leads to expressions in terms of the $m_{ij}$, as given by equation (8.9).

Thus we see that expressions for the difference $\bar{\pi}_j - \pi_j$ can be given in terms of the elements of generalized inverses $G = [g_{ij}]$, the elements of $H = [h_{ij}] = G(I - \Pi)$, the elements of the mean first passage times $M = [m_{ij}]$ and the elements of $N = [n_{ij}]$, The central role that the product $\pi_j m_{ij}$ plays, as reflected in the elements of $N$, leads to some interesting bounds (see Theorem 21 below.) The particular form that one wishes to use from Theorem 20 depends very much on the information that one has, be it a generalised inverse (of particular form), the mean first passage times (with and without the stationary probabilities) or the combination of the stationary probabilities and the mean first passage times through $n_{ij} = (1 - \delta_{ij})\pi_j m_{ij}$.

A thorough comparison of different bounds of the form $\|\bar{\pi}^T - \pi^T\| \leq \kappa_l \|E\|_q$ is given by Cho and Meyer (2001). Hunter, (2005) and Hunter, (2006) also explored the application of generalized inverse results to derive bounds on the stationary distributions under perturbations. We do not repeat that discussion here other than to make an observation regarding the following result given as Theorem 5.1 of Hunter, (2006):

**Theorem 21**: *For all irreducible m-state Markov chains undergoing a general perturbation* $E = [\varepsilon_{ij}]$,

$$\|\pi - \bar{\pi}\| = \sum_{j=1}^{m} |\pi_j - \bar{\pi}_j| \leq (K-1)\|E\|_\infty. \quad (8.6)$$

□

At the time of its derivation, the constant appearing in the bound was not identified as Kemeny's constant. Theorem 21 gives an interesting application of the role that Kemeny's constant provides in its involvement in overall bounds on the stationary probabilities subjected to perturbations. The reader is referred to the aforementioned papers for further results.

Using the results of the earlier sections of this paper we have been able to give some alternative derivations of earlier results and clarify their applicability.